\documentclass[a4paper,11pt]{article}
\usepackage[top=2.5cm,bottom=2.5cm,left=2.5cm,right=2.5cm]{geometry}
\usepackage{amsfonts}
\usepackage{latexsym,amsmath,amssymb,amsfonts,epsfig,psfrag,url,graphics,ifpdf,multicol}
\usepackage{color,colordvi,amscd,amsthm} 
\usepackage{tikz}
\usepackage[numbers,sort&compress]{natbib}
\usepackage{indentfirst,graphics,epsfig}
\usepackage{graphicx}
\usepackage{graphics}
\usepackage{graphicx,psfrag}
\usepackage{graphics}

\usepackage{tabularx}
\usepackage{booktabs}
\usepackage{floatrow}
\usepackage{multirow}
\usepackage{graphicx}
\usepackage{caption}
\usepackage[ruled]{algorithm2e}
\newfloatcommand{capbtabbox}{table}[][\FBwidth]

\newtheorem{theo}{Theorem}[section]

\newtheorem{coro}[theo]{Corollary}
\newtheorem{prob}[theo]{Problem}
 \theoremstyle{definition}

\theoremstyle{remark}

\newcounter{casenum}[theo]

\newcounter{subcasenum}[theo]

\newcounter{claimnum}[theo]

\allowdisplaybreaks

\usepackage{indentfirst,subfig}

\begin{document}
	\thispagestyle{empty}
	\captionsetup[figure]{labelfont={bf},name={Fig.},labelsep=period}

\begin{center} {\large\sc
A Survey on the $k$-Path Vertex Cover Problem}
\end{center}
\pagestyle{empty}
\begin{center}
{
  {\small Jianhua Tu$^{1,2,3}$}}\\[2mm]

{\small $^1$ School of Mathematics and Statistics, Beijing Technology and Business University, \\
		\hspace*{1pt} Beijing, P.R. China 100048} \\
{\small $^2$ Key Laboratory of Tibetan Information Processing and Machine Translation,\\
		\hspace*{1pt} Qinghai Province, XiNing, P.R. China 810008}\\
{\small $^3$ Key Laboratory of Tibetan Information Processing, Ministry of Education, \\
		\hspace*{1pt} XiNing, P.R. China, 810008}\\[2mm]
{	E-mail: tujh81@163.com}\\

\end{center}

\begin{center}

\begin{abstract}

Given a graph $G=(V,E)$ and a positive integer $k\ge2$, a $k$-path vertex cover is a subset of vertices $F$ such that every path on $k$ vertices in $G$ contains at least one vertex from $F$. A minimum $k$-path vertex cover in $G$ is a $k$-path vertex cover with minimum cardinality and its cardinality is called the {\it $k$-path vertex cover number} of $G$.
In the {\it $k$-path vertex cover problem}, it is required to find a minimum $k$-path vertex cover in a given graph. In this paper, we present a brief survey of the current state of the art in the study of the $k$-path vertex cover problem and the $k$-path vertex cover number.

\vspace{5mm}

\end{abstract}
\end{center}

\baselineskip=0.24in

\section{Introduction}

We consider only simple and undirected graphs and use standard graph theory notations (see
e.g. \cite{Bondy2008}). Let $G=(V,E)$ be a simple graph and $k\geq2$ be a positive integer. The number of vertices and edges in $G$ are called the order and size of $G$, respectively. For any subset $U\subseteq V$, let $G[U]$ denote the subgraph of $G$ induced by $U$, i.e.,  the subgraph of $G$ whose vertex set is $U$ and whose
edge set consists of all edges of $G$ which have both ends in $U$.
A path on $k$ vertices is called a $k$-path, or a $P_k$. Let $C_k$ and $K_k$ denote a cycle and a complete graph on $k$ vertices, respectively. Unless otherwise stated, we use $n$ and $m$ to refer to the number of vertices and edges, respectively, of the given graph.

A {\it $k$-path vertex cover} of in a graph $G$ is a subset of vertices $F$ such that every $k$-path in $G$ contains at least one vertex from $F$, that is, $G-F$ 
does not contain a path (not necessarily induced) on $k$ vertices. A $k$-path vertex cover is {\it minimum} if it has the smallest possible number of vertices. The {\it $k$-path vertex cover number}, $\psi_k(G)$, of $G$ is the cardinality of a minimum $k$-path vertex cover in $G$. In the {\it $k$-path vertex cover problem} (MinVCP$_k$), it is required to find a minimum $k$-path vertex cover in a given graph. In the literature, a $k$-path vertex is also called a vertex cover $P_k$ \cite{Tu2011-1,Tu2011-2}, or a vertex $k$-path cover \cite{Bresar2013}, or a $k$-observer \cite{Acharya2012,Ries2015}, or a $k$-path transversal \cite{Lee2019}, or a $P_k$-hitting set \cite{Brustle2021}.


Clearly, MinVCP$_2$ is exactly the well-known vertex cover problem that is a problem of finding a set of vertices with minimum cardinality in a graph that includes at least one endpoint of every edge of the graph. The vertex cover problem is one of the most important and fundamental computational problems in graph theory, combinatorial optimization, and theoretical computer science. It follows that $k$-path vertex cover is a natural generalization of the intensively studied concept of vertex cover. 

On the other hand, MinVCP$_k$ is also a special case of the vertex deletion problem \cite{Fujito1998,Krishnamoorthy1979, Yannakakis1981}. A graph property $\Pi$ is hereditary if it is closed under induced subgraphs, and non-trivial if it is satisfied by infinitely many graphs and it is not satisfied by infinitely many graphs. Given a hereditary and non-trivial graph property $\Pi$, the objective of the $\Pi$-vertex deletion problem is to find a vertex set with minimum cardinality whose deletion results in a graph satisfying $\Pi$. Many classic optimization problems like the feedback vertex set problem ($\Pi$ means ``containing no cycle"), the vertex bipartization problem ($\Pi$ means ``containing no odd cycle"), the $d$-bounded-degree vertex deletion problem ($\Pi$ means ``having maximum degree at most $d$") are examples of vertex deletion problems. MinVCP$_k$ corresponds to the case that $\Pi$ means ``containing no $k$-path".


Note that a graph that contains no 3-path has maximum degree at most 1, it follows that MinVCP$_2$ and MinVCP$_3$ are equivalent to the 0-bounded-degree and 1-bounded-degree vertex deletion problem, respectively. It is worth mentioning that a subset of vertices $F$ of a graph $G$ is a 3-path vertex cover if and only if its complement $V(G)\setminus F$ is a so-called dissociation set, that is, a subset of vertices inducing a subgraph with maximum degree at most 1. A maximum dissociation set in a graph $G$ is a dissociation set with maximum cardinality and its cardinality, denoted by $diss(G)$, is called the dissociation number of $G$. Clearly, $\psi_3(G)=|V(G)|-diss(G)$. The maximum dissociation set problem that is the dual problem of MinVCP$_3$ asks to find a maximum dissociation set in a given graph. The problem was introduced more than four decades ago by Yannakakis \cite{Yannakakis1981} and can be viewed as a generalization both of the maximum independent set problem and the maximum induced matching problem. The problem has been proved to be NP-hard in the class of bipartite graphs \cite{Yannakakis1981}.


The study of MinVCP$_k$ is motivated by the two real world problems that are related to secure communication in wireless sensor networks (WSNs) \cite{Bresar2011, Novotny2010} and to installation of traffic cameras \cite{Tu2011-1}, respectively.

Nowdays, WSNs have been widely used in various areas, such as
environment monitoring, disaster detection, and home automation, etc.
A WSN can be described by a graph composed of vertices and edges denoting, respectively, sensor devices and communication channels between pairs of sensor devices.
It is very important to design security protocols of WSNs in order to ensure security properties such as confidentiality, authenticity, and data integrity. The Canvas scheme that can provide data integrity in a sensor network was developed by Vogt \cite{Vogt2004}.
Novotn\'{y}\cite{Novotny2010} proposed the $k$-generalized Canvas scheme that can guarantee data integrity under the assumption that at least one node of each $k$-path is a protected node. Since a protected node is much more costly than an ordinary node, we need to minimize
the number of protected nodes. The problem of minimizing the number of protected nodes is exactly MinVCP$_k$.

Another motivation for studying MinVCP$_k$ is related to installation of traffic cameras. The increasing cars and buses lead to an increase in road traffic accidents, hence
it is very necessary to install cameras at traffic intersections. A road network can be modeled by a graph in which vertices represent traffic intersections and edges represent connections between pairs of traffic intersections. If every traffic intersection is installed with several cameras, the cost
would be enormous. For a given integer $k$, we aim to choose as few traffic intersections as possible to install cameras such that a driver will encounter at least one camera within $k$ traffic intersections. The corresponding optimization problem can be formulated as MinVCP$_k$. 

MinVCP$_k$ also finds applications in monitoring message flows in WSNs. Suppose every message which continuously passes $k$ nodes should be monitored at least once, then we again have MinVCP$_k$.

Basing on the real world problems mentioned above, three variants of MinVCP$_k$ have also been raised and studied.

\begin{itemize}
	\item {\it The weighted $k$-path vertex cover problem} ($\mathcal{W}$-MinVCP$_k$). Given a graph $G=(V, E)$ and a positive weight $w(v)$ for every vertex $v\in V$, the goal of $\mathcal{W}$-MinVCP$_k$ is to find a $k$-path vertex cover in $G$ with minimum total weight. 
	
	\item  {\it The connected $k$-path vertex cover problem} ($\mathcal{C}$-MinVCP$_k$). Given a connected graph $G=(V, E)$, the goal of $\mathcal{C}$-MinVCP$_k$ is to find a subset of vertices $F$ in $G$ with minimum cardinality such that $F$ is a $k$-path vertex cover of $G$ and the induced subgraph $G[F]$ is connected.
	
	\item  {\it The weighted connected $k$-path vertex cover problem} ($\mathcal{WC}$-MinVCP$_k$). Given a connected graph $G=(V, E)$ and a positive weight $w(v)$ for every vertex $v\in V$, the goal of $\mathcal{WC}$-MinVCP$_k$ is to find a subset of vertices $F$ in $G$ with minimum total weight such that $F$ is a $k$-path vertex cover of $G$ and the induced subgraph $G[F]$ is connected.
\end{itemize}

Due to their importance in theory and application, MinVCP$_k$ and its variants have been studied extensively. In particular, a large number of results on exact algorithms, approximation algorithms, and parameterized algorithms for MinVCP$_k$ and its variants have been reported. The present paper aims to provide a brief survey of the current state of the art in the study of MinVCP$_k$ and its variants, and the $k$-path vertex cover number. We mainly focus on the cases with $k\geq3$.


\section{Computational complexity}

The decision version of MinVCP$_k$ is stated as follows:

{\bf INPUT:} A graph $G$ and a positive integer $t$.

{\bf OUTPUT:} Is there a $k$-path vertex cover $F$ in $G$ of size at most $t$?

We abuse notation and let MinVCP$_k$ refer to the $k$-path vertex cover problem and its decision version.

\begin{theo}\cite{Bresar2011,Acharya2012}
 For any fixed integer $k\geq 2$, MinVCP$_k$ is NP-complete.
\end{theo}

The 2-subdivision of a graph $G$ is obtained from $G$ by 2-subdividing each edge of $G$, where to
2-subdivide an edge $e=uv$ is to delete $e$, add two new vertices $x$, $y$ and three news edges $ux$, $xy$, $yv$.
The 2-subdivision of a graph is called a 2-subdivision graph.
Poljak \cite{Poljak1974} showed that the vertex cover problem is NP-complete for 2-subdivision graphs. Following Poljak's result and the fact that 2-subdivision graphs are $C_t$-free, where $3\leq t\leq 8$, Brešar et al. \cite{Bresar2019} proved the following result.

\begin{theo}\cite{Bresar2019}
	For any fixed integer $k\geq 2$, MinVCP$_k$ is NP-complete for $C_t$-free graphs, where $3\leq t\leq 8$.
\end{theo}

The NP-completeness of MinVCP$_3$ has been studied intensively. Yannakakis \cite{Yannakakis1981} proved that MinVCP$_3$ remains NP-complete even in bipartite graphs. The author and Yang \cite{Tu2013} studied the NP-completeness of MinVCP$_3$ for cubic planar graphs. The girth of a graph is the length of one of its (if any) shortest cycles. Acyclic graphs are considered to have infinite girth.

\begin{theo}\cite{Tu2013}
MinVCP$_3$ is NP-complete for cubic planar graphs of girth 3.
\end{theo}

If P$\neq$NP is assumed, a wide class of optimization problems cannot be solved exactly in polynomial time and approximation algorithms naturally arise in the field of computer science and operations research. An {\it $\alpha$-approximation algorithm} for an optimization problem runs in polynomial time and outputs a solution that is within a factor $\alpha$ of being optimal. We call $\alpha$ the {\it approximation ratio} or {\it approximation factor} of the algorithm.

A {\it polynomial time approximation scheme} (PTAS) for an optimization problem is a set of algorithms such that, given $\varepsilon>0$, it takes an instance of the problem and produces a solution that is within a factor $1+\varepsilon$ of being optimal (or $1-\varepsilon$ for maximization problems) with running time polynomial in the instance size when $\varepsilon$ is fixed.
The class APX is the set of NP optimization problems that allow constant-factor approximation algorithms.
A problem is said to be {\it APX-hard} if there is a PTAS reduction from every problem in APX to that problem, and to be {\it APX-complete} if the problem is APX-hard and also in APX. If P$\neq$NP is assumed, no APX-hard problem has a PTAS.

\begin{theo}\cite{Kumar2014}
MinVCP$_3$ is APX-complete for bipartite graphs.
\end{theo}

\begin{theo}\cite{Devi2015}
MinVCP$_4$ is NP-complete for cubic planar graphs and APX-complete for cubic graphs, cubic bipartite graphs and $K_{1,4}$-free graphs.
\end{theo}

\section{Exact algorithms}

MinVCP$_k$ for trees can be solved in linear time.

\begin{theo}\cite{Bresar2011}
Let $T$ be a tree and $k$ be a positive integer. There exists a linear time algorithm that can return an optimal $k$-path vertex cover of $T$ of size at most $\frac{|V(T)|}{k}$. Therefore, $\psi_k(T)\leq\frac{|V(T)|}{k}$.
\end{theo}

Bre\v{s}ar et al. \cite{Bresar2014} considered $\mathcal{W}$-MinVCP$_k$, and gave linear time algorithms for $\mathcal{W}$-MinVCP$_k$ on complete graphs, cycles and trees.

\begin{theo}\cite{Bresar2014}
$\mathcal{W}$-MinVCP$_k$ on complete graphs $K_n$, cycles $C_n$ and trees with $n$ vertices can be solved with time complexity
$O(n\cdot k)$, $O(n\cdot k^2)$ and $O(n\cdot k)$, respectively.
\end{theo}

The Cartesian product $G=G_1 \square G_2$ of graphs $G_1$ and $G_2$ with disjoint vertex sets $V_1$ and $V_2$ and edge sets $E_1$ and $E_2$ is the graph with vertex set $V_1\times V_2$ and $u=(u_1,u_2)$ adjacent with $v=(v_1,v_2)$ whenever $u_1=v_1$ and $u_2v_2\in E_2$ or $u_2=v_2$ and $u_1v_1\in E_1$. A Cartesian product graph $P_m \square P_n$ of paths on $m$ and $n$ vertices is called a grid graph. MinVCP$_k$ for grid graphs can be solved in linear time \cite{Acharya2012}.

A cactus is a connected graph in which any two simple cycles have at most one vertex in common. Equivalently, every edge in such a connected graph belongs to at most one simple cycle. The author \cite{Tu2017} considered MinVCP$_k$ on cacti.

\begin{theo}\cite{Tu2017}
There exists an efficient algorithm that can obtain the $k$-path vertex cover number of a cactus $G$ in $O(n^2)$ time.
\end{theo}

Li et al. \cite{Li2016} considered $\mathcal{WC}$-MinVCP$_k$ and showed that $\mathcal{WC}$-MinVCP$_k$ can be solved in $O(n)$ time when the graph is a tree, and can be solved in $O(r\cdot n)$ time when the graph is a uni-cyclic graph, where $r$ is the length of the unique cycle.

There have been a lot of studies on exact algorithms for MinVCP$_2$ and MinVCP$_3$. MinVCP$_2$ can be solved in $\mathcal{O}^*(1.1996^n)$ time \cite{Xiao2017-4}. Throughout this paper, the $\mathcal{O}^*()$ notation suppresses all factors polynomial in the input size. Since the maximum dissociation set problem is the dual problem of MinVCP$_3$, in terms of exact algorithms, there's no need to distinguish these two problems. Exact algorithm for MinVCP$_3$ and the maximum dissociation problem has been improved for several times from the first 
$\mathcal{O}^*(1.5171^n)$-time algorithm due to Kardo\v{s} et al. \cite{Kardos2011}, to an $\mathcal{O}^*(1.4656^n)$-time algorithm due to Xiao and Kou \cite{Xiao2017-2}, and to an $\mathcal{O}^*(1.4613^n)$-time algorithm due to Chang et al. \cite{Chang2014, Chang2018}. Xiao and Kou \cite{Xiao2017-2} reduced the time complexity at the cost of an exponential space complexity and gave the following result.


\begin{theo}\cite{Xiao2017-2}
MinVCP$_3$ can be solved in $\mathcal{O}^*(1.3659^n)$ time and space.
\end{theo}

Given a universe set of $n$ elements and a collection $\mathcal{C}$ of subsets of size at most $k$, {\it the $k$-hitting set problem} asks for finding a hitting set of $\mathcal{C}$ with minimum cardinality, where a hitting set of $\mathcal{C}$ is a subset $U_0\subseteq U$ such that every subset of $\mathcal{C}$ contains at least one element of $U_0$. It is easy to see that MinVCP$_k$ is a special case of the $k$-hitting set problem. Following the results on the $k$-hitting set problem due to Fomin et al. \cite{Fomin2010} and Fernau \cite{Fernau2010}, one can obtain exact algorithms for MinVCP$_4$ and $\mathcal{W}$-MinVCP$_4$.

\begin{theo}\cite{Fomin2010}
MinVCP$_4$ can be solved in $\mathcal{O}^*(1.8704^n)$ time.
\end{theo}

\begin{theo}\cite{Fernau2010}
	$\mathcal{W}$-MinVCP$_4$ can be solved in $\mathcal{O}^*(1.97^n)$ time.
\end{theo}

Many special classes of graphs are known where MinVCP$_3$ is solvable in polynomial time. Together with the complexity of the maximum dissociation set problem, it is known that MinVCP$_3$ is polynomially solvable for chordal and weakly chordal graphs, asteroidal triple-free (AT-free) graphs \cite{Cameron2006}, $(P_k, K_{1,n})$-free graphs (for any positive $k$ and $n$) \cite{Lozin2003}, (claw, bull)-free graphs, (chair, bull)-free graphs, (chair, $K_3$)-free graphs, $\ell K_2$-free graphs (for any fixed integer $\ell\geq2$) \cite{Orlovich2011}, $P_5$-free graphs \cite{Brause2017} and some other hereditary classes of graphs \cite{Cameron2006, Lozin2003, Boliac2004}.
In particular, MinVCP$_3$ can be solved in linear time for $P_4$-tidy graphs \cite{Bresar2019} and line graphs of graphs having a
Hamiltonian path \cite{Orlovich2011}.

\section{Approximation algorithms}

A trivial $k$-approximation algorithm for MinVCP$_k$ can be easily obtained. One way to find a $k$-path vertex cover is to repeat the following process: find a $k$-path, put its vertices into solution $S$, and remove all edges incident to any vertex in $S$ from the graph. As any $k$-path vertex cover of the input graph must contain at least one vertex of each $k$-path that was considered in the process, the solution produced, therefore, is within a factor $k$ of the optimal one.  


Bre\v{s}ar et al. \cite{Bresar2011} proved that for any $\alpha\geq1$, an $\alpha$-approximation algorithm for MinVCP$_k$, with
polynomial running time, yields directly an $\alpha$-approximation for MinVCP$_2$. MinVCP$_2$ is APX-complete even on cubic graphs \cite{Alimonti2000}. Furthermore, MinVCP$_2$ cannot be approximated with $\sqrt{2}$ unless P=NP \cite{Khot2002, Khot2018} and cannot be approximated within any constant factor less than 2 under the unique game conjecture \cite{Khot2002, Khot2008}. Thus for every $k\geq3$, it is NP-hard to approximate MinVCP$_k$ within $\sqrt{2}$ unless P=NP.

Ries et al. \cite{Ries2015} presented a 3-approximation algorithm for MinVCP$_k$ on $d$-regular graphs for $k\leq \frac{d+2}{2}$. Their algorithm is the first one with an better approximation factor than $k$ for general $k$. Zhang et al. \cite{Zhang2020} gave an improved result.

\begin{theo}\cite{Zhang2020}
	When $1\leq k-2 < d$, there is a $\frac{\lfloor d/2\rfloor(2d-k+2)}{(\lfloor d/2\rfloor+1)(d-k+2)}$-approximation algorithm for MinVCP$_k$ on $d$-regular graphs that runs in $O(d^2\cdot n)$ time.
\end{theo}

Note that when $k\leq \frac{d+2}{2}$, $\frac{\lfloor d/2\rfloor(2d-k+2)}{(\lfloor d/2\rfloor+1)(d-k+2)}<3$. Lee \cite{Lee2019} considered MinVCP$_k$ on general graphs and presented an approximation algorithm that strictly improves the trivial $k$-approximation algorithm for general graphs.

\begin{theo}\cite{Lee2019}
There is an $O(\log k)$-approximation algorithm for MinVCP$_k$ that runs in $2^{O(k^3 \log k)}n^2 \log n+ n^{O(1)}$ time.
\end{theo}

When $k$ is a constant, the algorithm runs in polynomial-time. Note that for any approximation algorithm, the running time cannot be polynomial in $k$ since detecting a single $k$-path is NP-hard (it includes Hamiltonian Path as a special case).  

Very recently, Br\"{u}stle et al. \cite{Brustle2021} showed that MinVCP$_k$ admits a $(k-\frac{1}{2})$-approximation. In fact, they studied a wider problem called the $H$-hitting set problem. Let $H$ be a fixed graph of order $k$, the goal of the $H$-hitting set problem is to delete a minimum number of vertices in a given graph $G$ such that the resulting graph has no copies of $H$ as a subgraph. Clearly, the $P_k$-hitting set problem is exactly MinVCP$_k$. Br\"{u}stle et al. \cite{Brustle2021} proved that for a tree $T$ with $k$ vertices, the $T$-hitting set problem admits a $(k-\frac{1}{2})$-approximation. Thus, we have

\begin{theo}\cite{Brustle2021}
MinVCP$_k$ admits a $(k-\frac{1}{2})$-approximation.
\end{theo}

A graph $G$ is a $d$-dimensional ball graph if each vertex of $G$ corresponds to a ball in $\mathbb{R}^d$, two vertices are adjacent in $G$ if and only if their corresponding balls have nonempty intersection. Let $r_{max}$ and $r_{min}$ be the largest and smallest radius of those balls, respectively. We write $r_{max}/r_{min}$ to denote the heterogeneity of a ball graph. A 2-dimensional ball graph is also known as a disk graph. A unit disk graph is a disk graph in which all disks have the same radii. 
For fixed integer $k$, Zhang et al. \cite{Zhang2017} present a PTAS for MinVCP$_k$ on a ball graph whose heterogeneity is bounded by a constant.

Liu et al. \cite{Liu2013} ] were the first to study $\mathcal{C}$-MinVCP$_k$ and presented a PTAS for $\mathcal{C}$-MinVCP$_k$ on unit disk graphs using partition technique and shifting strategy. A simpler PTAS for $\mathcal{C}$-MinVCP$_k$ on unit disk graphs was given in \cite{Chen2018}. The simpler PTAS not only simplifies the previous algorithm, but also reduces the time-complexity. 

Li et al. \cite{Li2016} showed that $\mathcal{C}$-MinVCP$_k$ can be approximable within $k$ assuming that the graph has girth at least $k$. Fujito \cite{Fujito2017} improved this result and showed that when $k$ is a fixed parameter, $\mathcal{C}$-MinVCP$_k$ admits a $k$-approximation for general graphs without any assumption on girth. He also studied $\mathcal{WC}$-MinVCP$_k$. 

\begin{theo}\cite{Fujito2017}
For any fixed integer $k\geq2$, $\mathcal{C}$-MinVCP$_k$ admits a $k$-approximation, and $\mathcal{WC}$-MinVCP$_k$ is as hard to approximate as the weighted set cover problem, but approximable within $1.35\ln n + 3$ for $k\leq 3$.
\end{theo}


In the last decade, there has been a great deal of interest in the study of approximation algorithms for MinVCP$_3$, MinVCP$_4$ and their variants. Kardo\v{s} et al. \cite{Kardos2011} gave a randomized algorithm for MinVCP$_3$ with expected approximation ratio $\frac{23}{11}$. The author and Zhou presented two 2-approximation algorithms for $\mathcal{W}$-MinVCP$_3$ using primal–dual method \cite{Tu2011-1} and using local ratio method \cite{Tu2011-2}.

\begin{theo}\cite{Tu2011-1,Tu2011-2}
	There exist 2-approximation algorithms with running time $O(m\cdot n)$ for $\mathcal{W}$-MinVCP$_3$.
\end{theo}

Chang et al. \cite{Chang2018} gave a moderately exponential time approximation algorithm for MinVCP$_3$.

\begin{theo}\cite{Chang2018}
There exists a $\frac{4}{3}$-approximation algorithm for MinVCP$_3$ that runs in $\mathcal{O}^*(1.4159^n)$ time.
\end{theo} 

Camby et al. \cite{Camby2014} gave a 3-approximation algorithm for MinVCP$_4$ using primal-dual method. On the other hand, there are a lot of works on approximation algorithms for MinVCP$_3$ and MinVCP$_4$ on some special classes of graphs \cite{Tu2013,Li2014,Devi2015,Ries2015,Zhang2020}. Table \ref{table1} summarizes the best approximation ratios for MinVCP$_3$ and MinVCP$_4$.

\begin{center}
	\begin{table}[h]
		\caption{The best approximation ratios for MinVCP$_3$ and MinVCP$_4$}\label{table1}
		\renewcommand\arraystretch{2}
		\footnotesize
		\begin{tabular}{|m{12mm}<{\centering}|m{12mm}<{\centering}m{12mm}<{\centering}m{15mm}<{\centering}m{14mm}<{\centering}m{18mm}<{\centering}m{24mm}<{\centering}m{13mm}<{\centering}|}
			\hline
			 MinVCP$_k$&General graphs&Bipartite graphs&Cubic graphs&4-regular graphs&Bipartite $d$-regular graphs&$d$-regular graphs  $(d\geq5)$&$K_{1,4}$-free graphs\\
			\hline
			$k=3$&2 \cite{Tu2011-1,Tu2011-2}&&$1.25$ \cite{Zhang2020}&$\frac{14}{9}$ \cite{Zhang2020}&$\frac{2d-1}{2d-2}$ \cite{Ries2015}&$\frac{\lfloor\frac{d}{2}\rfloor(2d-1)}{(\lfloor\frac{d}{2}\rfloor+1)(d-1)}$ \cite{Zhang2020}&\\
		\hline
			$k=4$&3 \cite{Camby2014}&2 \cite{Kumar2014}&$\frac{15}{8}+\varepsilon$ for any $\varepsilon>0$ \cite{Zhang2020}&1.852 \cite{Zhang2020}&$\frac{d^2}{d^2-d+1}$ \cite{Zhang2020}&$\frac{(3d-2)(2d-2)}{(3d+4)(d-2)}$ (when $d$ is even) \cite{Zhang2020}&3 \cite{Devi2015}\\
			\hline
		\end{tabular}
	\end{table}
\end{center}

An {\it efficient polynomial time approximation scheme} (EPTAS) for an optimization problem $Q$ is a PTAS such that, given an instance $I$ of $Q$ and $\varepsilon>0$, its runtime is bounded by $O(f(\varepsilon)|I|^c)$, where $f$ is an arbitrary function and $c$ is a constant.
Using the dynamic programming algorithm and the Baker’s EPTAS framework for NP-hard
problems \cite{Baker1994}, the author and Shi \cite{Tu2019} presented an EPTAS for MinVCP$_3$ on planar graphs.

For $\mathcal{C}$-MinVCP$_3$, the problem admits a 3-approximation \cite{Li2016,Fujito2017}. Liu et al. \cite{Liu2020} presented a polynomial-time $(2\alpha+1/2)$-approximation algorithm for $\mathcal{C}$-MinVCP$_3$, where $\alpha$ is the approximation factor of an approximation algorithm for MinVCP$_3$. For those classes of graphs on which MinVCP$_3$ can be approximable within $\alpha<5/4$, their algorithm is a kind of improvement.

For $\mathcal{WC}$-MinVCP$_3$, Ran et al. \cite{Ran2019} presented a $(\ln \Delta+4+\ln 2)$-approximation algorithm, where $\Delta$ is the maximum degree of the input graph.  
Wang et al. \cite{Wang2015} showed that the problem remains NP-hard when restricted
to unit disk graphs and gave a PTAS for the problem on unit disk graphs under the assumption that the problem is $c$-local and the unit disk graphs have minimum degree of at least two. Wang et al. \cite{Wang2017} gave a PTAS for the problem on unit ball graphs under the assumption that the weight is smooth and weak $c$-local.

\section{Parameterized algorithms}

Parameterized algorithmics analyzes running time in finer detail than classical complexity theory: instead of expressing the running time as function of the input size only, dependence on one or more parameters of the input instance is taken into account. 

In parameterized algorithmics, the parameter $t$ may be the size of the solution sought after, or a number describing how ``structured'' the input instance is. 
A {\it fixed-parameter algorithm} (parameterized algorithm, or FPT algorithm) of a parameterized problem $(I,t)$ is an exact algorithm with running time $f(t)\cdot|I|^c$, for a constant $c$ independent of both $|I|$ and $t$. Thus, fixed-parameter algorithms confine the exponential part of the running time to the parameter $t$, presumably much smaller than the input size. A parameterized problem that admits a fixed-parameter algorithm is called {\it fixed parameter tractable} with respect to the parameter $t$.

When parameterized by the size $t$ of the solution sought after, there exists a trivial FPT algorithm for MinVCP$_k$ that runs in $\mathcal{O}^*(k^t)$ time \cite{Cai1996}. As mentioned above, MinVCP$_k$ is a special case of the $k$-hitting set problem. According to the results on the $k$-hitting set problem due to Fomin et al. \cite{Fomin2010}, the following result holds. 

\begin{theo}\cite{Fomin2010}
For any $k\geq 4$, MinVCP$_k$ can be solved in time $\mathcal{O}^*((k-0.9245)^t)$. 
\end{theo}

For MinVCP$_2$, the current best known algorithm due to Chen et al. \cite{Chen2010} runs in time $\mathcal{O}^*(1.2738^t)$. 
An $\mathcal{O}^*(2^t)$-time FPT algorithm for MinVCP$_3$ was given by the author \cite{Tu2015}. This result was subsequently improved for several times from an $\mathcal{O}^*(1.882^t)$-time algorithm due to Wu \cite{Wu2015}, to an $\mathcal{O}^*(1.8172^t)$-time algorithm due to Katreni\v{c} \cite{Katrenic2016}, to an algorithm due to Chang et al. \cite{Chang2016} running in $\mathcal{O}^*(1.7485^t)$ time and in exponential space, to an algorithm due to Xiao and Kou \cite{Xiao2017-3} running in $\mathcal{O}^*(1.7485^t)$ time and in polynomial space, and to an $\mathcal{O}^*(1.713^t)$-time algorithm due to Tsur \cite{Tsur2019}.


The author and Jin \cite{Tu2016} gave an $\mathcal{O}^*(3^t)$-time FPT algorithm for MinVCP$_4$ using iterative compression method. Tsur \cite{Tsur2021} improved the result to $\mathcal{O}^*(2.619^t)$.
For $k=5,\ 6$, and $7$, FPT algorithms for MinVCP$_k$ were also given in \cite{Cervenvy2019, Tsur2019-2}. Very recently, \v{C}erven\'{y} and Such\'{y} gave FPT algorithms outperforming those previously known for MinVCP$_k$ for $3\leq k \leq 8$ by developing a framework to automatically generate parameterized branching algorithms. The improved results given in \cite{Cervenvy2021} are summarized in Table \ref{table3}.

\begin{center}
	\begin{table}[h]
		\caption{Improved running times of MinVCP$_k$ for $3\leq k \leq 8$}\label{table3}
		\renewcommand\arraystretch{1.5}
		\footnotesize
		\begin{tabular}{|m{15mm}<{\centering}|m{15mm}<{\centering}m{15mm}<{\centering}m{15mm}<{\centering}m{15mm}<{\centering}m{15mm}<{\centering}m{15mm}<{\centering}|}
			\hline
			MinVCP$_k$&$k=3$&$k=4$&$k=5$&$k=6$&$k=7$&$k=8$\\
			\hline
			&$\mathcal{O}^*(1.712^t)$&$\mathcal{O}^*(2.151^t)$&$\mathcal{O}^*(2.695^t)$&$\mathcal{O}^*(3.45^t)$&$\mathcal{O}^*(4.872^t)$&$\mathcal{O}^*(5.833^t)$\\
			\hline
		\end{tabular}
	\end{table}
\end{center}

For MinVCP$_3$ on planar graphs, the author et al. \cite{Tu2017-2} showed that there is a subexponential parameterized algorithm for the problem running in time $\mathcal{O}^*(2^{O(\sqrt{t})})$.

The author et al. \cite{Tu2017-2, Bai2019} also proved that MinVCP$_3$ is fixed parameter tractable with respect to the treewidth of the input graphs. The treewidth is a fundamental graph parameter which captures how similar a graph is to a tree. Given an $n$-vertex graph together with its tree decomposition of width at most $p$, an algorithm running in time $\mathcal{O}^*(3^p)$ for MinVCP$_3$ was given in \cite{Bai2019}.

Consider $\mathcal{W}$-MinVCP$_k$. Shachnai and Zehavi \cite{Shachnai2015} introduced a multivariate approach for solving weighted
parameterized problems. Following their results, one can obtain FPT algorithms for $\mathcal{W}$-MinVCP$_2$ and $\mathcal{W}$-MinVCP$_3$.

\begin{theo}\cite{Shachnai2015}
For a graph $G=(V,E)$, a weight function $w: V\rightarrow [1,+\infty)$, and a parameter $W\geq 1$, there exists an FPT algorithm for $\mathcal{W}$-MinVCP$_2$ that runs in $\mathcal{O}^*(1.381^s)$ time and in polynomial space, or in $\mathcal{O}^*(1.363^s)$ time and space, where $s\leq W$ is the minimum size of a solution of
weight at most $W$.
\end{theo}

\begin{theo}\cite{Shachnai2015}
	For a graph $G=(V,E)$, a weight function $w: V\rightarrow [1,+\infty)$, and a parameter $W\geq 1$, there exists an FPT algorithm for $\mathcal{W}$-MinVCP$_3$ that runs in $\mathcal{O}^*(2.168^s)$ time and in polynomial space, where $s\leq W$ is the minimum size of a solution of
	weight at most $W$.
\end{theo}

Fernau \cite{Fernau2010} studied parameterized algorithms for the weighted $k$-hitting set problem. Following his results, one can obtain FPT algorithms for $\mathcal{W}$-MinVCP$_k$ for any $k\geq4$.


\begin{theo}\cite{Fernau2010}
	For a graph $G=(V,E)$ and a weight function $w: V\rightarrow [1,+\infty)$, $\mathcal{W}$-MinVCP$_4$ can be solved in $\mathcal{O}^*((3.1479)^W)$ time, where $W$ is the total weight of the solution sought after.
\end{theo}

\begin{theo}\cite{Fernau2010}
	For a graph $G=(V,E)$ and a weight function $w: V\rightarrow [1,+\infty)$, $\mathcal{W}$-MinVCP$_k$ can be solved in $\mathcal{O}^*((c_k)^W)$ time, where $c_k$ is the largest positive root of the characteristic polynomial \[x^4-3x^3-(k^2-5k+5)x^2 +x +(k^2-6k+9),\] and $W$ is the total weight of the solution sought after. Some values of $c_k$ are listed below:
	\begin{center}
		\begin{table}[h]
			\caption{Some values of $c_k$}\label{table2}
			\renewcommand\arraystretch{1.5}
			\footnotesize
			\begin{tabular}{|m{15mm}<{\centering}|m{10mm}<{\centering}m{10mm}<{\centering}m{10mm}<{\centering}m{10mm}<{\centering}m{10mm}<{\centering}m{10mm}<{\centering}m{10mm}<{\centering}|}
				\hline
				$k$&$5$&$6$&$7$&$8$&$9$&10&100\\
				\hline
				$c_k\ (\leq)$&$4.1017$&$5.0640$&$6.0439$&$7.0320$&$8.0243$&$9.0191$&$99.0002$\\
				\hline
			\end{tabular}
		\end{table}
	\end{center}
\end{theo}

Kernelization is an important method that is often used for dealing with NP-hard problems. We say that two instances of a decision
problem are equivalent if and only if they are either both yes-instances or both no-instances.
For a given parameterized problem $\Pi$ with a parameter $t$, a {\it kernelization algorithm} is a polynomial-time algorithm which, for an input instance $(I,t)$ of $\Pi$, returns an equivalent instance $(I',t')$, called a {\it kernel}, such that $t'\leq t$ and the size of $I'$ is bounded by a function $g(t)$. If $g(t)$ is a linear function we call the kernel a {\it linear kernel}. 

Since an instance of MinVCP$_k$ can be formulated as an instance of the $k$-hitting set problem, MinVCP$_k$ admits a kernel with $O(t^{k-1})$ vertices and $O(k\cdot t^k)$ edges \cite{Abu-Khzam2010, Fafianie2015}, where $t$ is the size of the solution sought after. On the other hand, it is not possible to achieve a kernel with $O(t^{2-\varepsilon})$ edges for MinVCP$_k$ unless coNP is in NP/poly \cite{Dell2010}. Very recently, \v{C}erven\'{y} et al. \cite{Cervenvy2021-2} gave a kernel with $O(t^3)$ edges for MinVCP$_k$ for any $k\geq 6$.

\begin{theo}\cite{Cervenvy2021-2}
For any $k\geq6$, MinVCP$_k$ admits a kernel with $O(t^3k^{O(k)})$ vertices and edges.
\end{theo}

For MinVCP$_2$, the current best kernel known is due to Lampis \cite{Lampis2011} with $2t-c\log t$ vertices for any fixed constant $c$. Nemhauser and Trotter \cite{Nemhauser1975} proved a famous theorem (the NT-Theorem) for MinVCP$_2$.

\begin{theo}\cite{Nemhauser1975}
For a graph $G=(V,E)$ with $n$ vertices and $m$ edges, there exists an $O(\sqrt{mn}+n)$-time
algorithm to compute two disjoint vertex subsets $A$, $B$, such that the following three properties hold:
	
(1) If $S'$ is a vertex cover of $G[V\setminus (A\cup B)]$, then $A\cup S'$ is a vertex cover of $G$.
	
(2) There is a minimum vertex cover $S$ of $G$ with $A\subseteq S.$

(3) Every vertex cover of $G[V\setminus (A\cup B)]$ has size at least $|V\setminus (A\cup B)|/2.$
\end{theo}

Fellows et al. \cite{Fellows2011} extended the NT-Theorem for the $d$-bounded-degree vertex deletion problem and Xiao \cite{Xiao2017-1} improved Fellow et al.'s result. Since MinVCP$_3$ is equivalent to the 1-bounded-degree vertex deletion problem, one can derive a generalization of the NT-theorem for MinVCP$_3$.

\begin{theo}\cite{Xiao2017-1}
For a graph $G=(V,E)$ with $n$ vertices and $m$ edges, there exists an $O(n^{5/2}m)$-time
algorithm to compute two disjoint vertex subsets $A$ and $B$, such that the following three properties hold:
	
(1) If $S'$ is a 3-path vertex cover of $G[V\setminus (A\cup B)]$, then $A\cup S'$ is a 3-path vertex cover of $G$.
	
(2) There is a minimum 3-path vertex cover $S$ of $G$ with $A\subseteq S.$
	
(3) Every 3-path vertex cover of $G[V\setminus (A\cup B)]$ has size at least $|V\setminus (A\cup B)|/13.$
\end{theo}

The generalization of the NT-theorem for MinVCP$_3$ implies a kernel with $13t$ vertices for MinVCP$_3$. The bound of the size of the kernel for MinVCP$_3$ was subsequently improved for several times \cite{Brause2016, Xiao2017-3}. The current best kernel known is due to Xiao and Kou \cite{Xiao2017-3} with $5t$ vertices.

\v{C}erven\'{y} et al. \cite{Cervenvy2021-2} also gave kernels with $O(t^2)$ edges for MinVCP$_4$ and MinVCP$_5$ that are asymptotically optimal (unless coNP is in NP/poly).

\begin{theo}\cite{Cervenvy2021-2}
MinVCP$_4$ admits a kernel with at most $176t^2 + 166t$ edges.
\end{theo}

\begin{theo}\cite{Cervenvy2021-2}
	MinVCP$_5$ admits a kernel with  at most $608t^2+583t$ edges.
\end{theo}

\section{The $k$-path vertex cover number}

The decision version of the $k$-path vertex cover number $\psi_k(G)$ is NP-complete, moreover, in the case $k=3$ it is NP-complete even in $C_4$-free bipartite graphs with maximum degree 3 \cite{Boliac2004}. Some lower and upper bounds on $\psi_k(G)$ in terms of various graph parameters were given in the literature.

\begin{theo}\cite{Bresar2013}
Let $k\geq 2$ and $n\geq k$ be positive integers. Then $\psi_k(P_n)=\lfloor\frac{n}{k}\rfloor$, $\psi_k(C_n)=\lceil\frac{n}{k}\rceil$ and $\psi_k(K_n)=n-k+1$.
\end{theo}

\begin{theo}\cite{Li2018}
Let $n\geq 4$ be an integer. Then the following hold.

(i) $\psi_k(K_{n,n})=\lceil \frac{n+1}{2} \rceil$ for $k=n+1$;

(ii) $\psi_k(K_{n,n})=n+1-\lfloor \frac{k}{2} \rfloor$ for $n+2\leq k<2n$;

(iii)  $\psi_k(K_{n,n})\geq\frac{n^2-nk+2n}{n-\frac{k}{2}+1}$ for $4\leq k\leq n$.
\end{theo}
	
\begin{theo}\cite{Li2018}
If $m\geq n\geq2$, $\psi_2(K_{m,n})=\psi_3(K_{m,n})=n$. If $m>n\geq2$ and $k\geq 3$, we have
\[\psi_k(K_{m,n})=n+1-\lfloor \frac{k}{2}\rfloor.\]
\end{theo}
	
\begin{theo}\cite{Bresar2011}
Let $T$ be a tree. Then, $\psi_k(T)\leq |V(T)|/k$.
\end{theo}

\begin{theo}\cite{Bresar2013}
Let $k\geq 2$ and $d\geq k-1$ be positive integers. Then, for any $d$-regular graph $G$, the following holds:
\[\psi_k(G)\geq \frac{d-k+2}{2d-k+2}|V(G)|.\]
\end{theo}

Recently, Bujt\'{a}s et al. \cite{Bujtas2021} generalized the result to arbitrary graphs in terms of minimum and maximum degree.

\begin{theo}\cite{Bujtas2021}
Let $G$ be a graph with minimum degree $\delta$ and maximum degree $\Delta$. If $k\geq 3$ and $\delta\geq k-1$, then
	\[\psi_k(G)\geq \frac{\delta-k+2}{\delta+\Delta-k+2} |V(G)|.\]
\end{theo}

\begin{theo}\cite{Bujtas2021}
Let $k$ and $\Delta$ be integers with $k\geq3$ and $\Delta=2$ or $\Delta\geq 4$, and let $G$ be
a graph of maximum degree at most $\Delta$. Then the following hold.

(i) If $\Delta\geq 2$ is even, then \[\psi_k(G)\leq \frac{(k-1)(\Delta-2)+4}{(k-1)\Delta+4} |V(G)|.\]

(ii) If $\Delta\geq 5$ is odd, then \[\psi_k(G)\leq \frac{(k-1)(\Delta-3)+8}{(k-1)(\Delta-1)+8} |V(G)|.\]
\end{theo}

\begin{theo}\cite{Bujtas2021}
Let $G$ be a graph without isolated vertices. Then for every integer $k\geq 3$, the following holds.
	\[\psi_k(G)\leq |V(G)|-\frac{2k-3}{k-1}\sum_{v\in V(G)}\frac{1}{1+d(v)}.\]
\end{theo}

A graph is called chordal if it does not contain any induced cycles of order at least 4.

\begin{theo}\cite{Bujtas2021}
Let $G$ be a chordal graph with clique number $\omega$ and $k\geq 3$ be an integer. Then
\[\psi_k(G)\leq \frac{\omega}{\omega+k-1}|V(G)|.\]
\end{theo}
	
The lower and upper bounds on $\psi_k(G)$ for small $k$ were also investigated.

\begin{theo}\cite{Bresar2011, Bujtas2021}
Let $G$ be a graph with $n$ vertices and $m$ edges. Then 

(i) $\psi_3(G)\leq \max\{\frac{2n+m}{6},\frac{n+m}{4},\frac{4n+m}{9}\}$,

(ii) $\psi_4(G)\leq \frac{n+3m}{10}$.
\end{theo}

\begin{theo}\cite{Bresar2011}
Let $G$ be a graph of order $n$ and of maximum degree $\Delta$. Then
	\[\psi_3(G)\leq \frac{\lceil\frac{\Delta-1}{2}\rceil}{\lceil\frac{\Delta+1}{2}\rceil}n.\]
\end{theo}

\begin{coro}\cite{Tu2013}
Let $G$ be a cubic graph of order $n$, then $2n/5 \leq \psi_3(G)\leq n/2$.
\end{coro}

\begin{theo}\cite{Bresar2013}
Let $G$ be a graph of order $n$, size $m$ and average degree $d(G)$, and $\ell$ be the smallest positive integer such that
$d(G)\leq 2\ell + 2$. Then
\[\psi_3(G) \leq\frac{\ell n}{\ell+2}+ \frac{m}{(\ell+1)(\ell+2)}.\]
\end{theo}

\begin{theo}\cite{Bujtas2021}
Let $G$ be a planar graph of order $n$, then $\psi_3(G)\leq 11n/15$ and $\psi_6(G)\leq 2n/3.$
\end{theo}

\begin{theo}\cite{Bresar2011}
Let $G$ be an outerplanar graph of order $n$. Then $\psi_3(G)\leq n/2.$
\end{theo}

\begin{theo}\cite{Bujtas2021}
Let $G$ be a chordal and planar graph of order $n$. Then $\psi_3(G)\leq 2n/3$, $\psi_4(G)\leq 4n/7$, and $\psi_5(G)\leq n/2.$
\end{theo}

Some bounds and exact values for $\psi_k(G)$ on Cartesian product graphs, rooted product graphs, and Kneser graphs were also given \cite{Li2018, Bresar2013, Jakovac2015, Bresar2021}.

\section{Conclusions}

MinVCP$_k$ is a natural generalization of the classic vertex cover problem and has received increasing attention in recent years. This work aims to present
a survey on the $k$-path vertex cover problem that could be used to gain insights on the topic. In this survey, we mainly focus on the cases with $k\geq3$. It is important to point out that although there are various methods to deal with NP-hard problems, we focus mainly on the results on exact algorithms, approximation algorithms and parameterized algorithms for MinVCP$_k$ and its variants. In fact, some results on online algorithms and heuristic algorithms for MinVCP$_k$ have also been reported in the literature \cite{Zhang2019, Zhang2020-2}. Finally, we present some important problems which could be used to draw future research directions in this area.

\begin{prob}
Find $(1-\delta)k$-approximation for MinVCP$_k$ for a general $k$ and a universal constant $\delta>0$.
\end{prob}

\begin{prob}
Determine whether MinVCP$_k$ has an FPT algorithm running in $O^*((k-1-\varepsilon)^t)$ time for every value of $k$, where $t$ is the size of solution sought after.
\end{prob}

\begin{prob}
Determine whether the parameterized problem of MinVCP$_k$ has a kernel with a linear number of vertices, or a kernel with $O(t^2)$ edges for every value of $k$, where $t$ is the size of solution sought after.
\end{prob}

\begin{prob}
Find a generalization of the NT-theorem for MinVCP$_k$ for general $k$.
\end{prob}

\begin{prob}
 Prove and disprove the following conjecture: If $G$ is a planar graph of order $n$, then $\psi_3(G)\leq2n/3$.
\end{prob}

\section*{Acknowledgments}

The work was supported by Research Foundation for Advanced Talents of Beijing Technology and
Business University (No. 19008021187).

\bibliographystyle{unsrt}

\section*{Statements and Declarations}

\subsection*{Funding}

The work was supported by Research Foundation for Advanced Talents of Beijing Technology and
Business University (No. 19008021187).

\subsection*{Competing Interests}

The author has no relevant financial or non-financial interests to disclose.
\end{document}